\documentclass[12pt,reqno]{amsart}

\usepackage[arrow,matrix,curve]{xy}

\usepackage[dvips]{graphicx} 

\usepackage{amssymb, latexsym, amsmath, amscd, array, hyperref,
epigraph
}

\usepackage{breakurl}   

\usepackage[LGR,T1]{fontenc} 
\usepackage{wrapfig}  
\usepackage{lscape}   
\usepackage{rotating} 

\theoremstyle{definition}

\numberwithin{equation}{section}

\numberwithin{equation}{section}

\title{Three case studies in current Leibniz scholarship}

\author[M. Katz]{Mikhail G. Katz} \address{M.\;Katz, Department of
  Mathematics, Bar Ilan University, Ramat Gan 5290002 Israel}
\email{katzmik@math.biu.ac.il}

\author[K.~Kuhlemann]{Karl Kuhlemann}\address{K.\;Kuhlemann, Institute
  of Mathematics and Physics Education, Gottfried Wilhelm Leibniz
  University Hannover, D-30167 Hannover, Germany}
\email{kus.kuhlemann@t-online.de}

\author[D.~Sherry]{David Sherry} \address{D.\;Sherry, Department of
  Philosophy, Northern Arizona University, Flagstaff, AZ 86011, US}
\email{David.Sherry@nau.edu}

\author[M.~Ugaglia]{Monica Ugaglia} \address{M.\;Ugaglia, Il Gallo
  Silvestre, Localit\`a Collina 38, Montecassiano, Italy}
\email{monica.ugaglia@gmail.com}


\begin{document}

\thispagestyle{empty}



\begin{abstract}

We examine some recent scholarship on Leibniz's philosophy of the
infinitesimal calculus.  We indicate difficulties that arise in
articles by Bassler, Knobloch, and Arthur, due to a denial to
Leibniz's infinitesimals of the status of mathematical entities
violating Euclid V Definition 4.
\end{abstract}

\maketitle

\today

\tableofcontents

\epigraph{Infinitesimal calculus is as much, or as little, about
  in\-fi\-ni\-tesimals as the business of a variety shop is to sell
  varieties.  --Bassler}

\section{Three case studies}
\label{s6}

Leibniz's references to infinity and infinitesimals have been a source
of continued debate in current Leibniz scholarship; see e.g.,
\cite{22b} for a comparison of rival viewpoints.  Ishiguro developed a
so-called syncategorematic interpretation of Leibnizian infinities and
infinitesimals in her \cite[Chapter\;5]{Is90}.  A rival interpretation
was developed by Esquisabel and Raffo Quintana in \cite[pp.\;620,
  641]{Es21}.

The bone of contention is the status of Leibnizian infinitesimals: are
they mathematical entities?  Ishiguro denies this, and seeks to
interpret each occurrence of the term \emph{infinitesimal} in Leibniz
as merely shorthand for more long-winded exhaustion-style arguments
exploiting only ``ordinary'' numbers.  Such a state of affairs is
sometimes referred to as \emph{syncategorematicity}.  As in Mediaeval
logic \emph{syncategoremata} were the parts of speech lacking a
definite reference of their own and depending instead on other parts
(\emph{categoremata}) in order to acquire a definite reference,
syncategorematic infinitesimals depend on other mathematical entities,
namely the ones involved in the method of exhaustion.  See also
\cite[Section~3]{22a}.

Esquisabel and Raffo Quintana, while acknowledging the fictional
nature of Leibnizian infinitesimals, nonetheless accord to them the
status of mathematical entities.  Bascelli et al.\;\cite{16a},
B{\l}aszczyk et al.\;\cite{17e}, and Katz et al.\;\cite{22a} similarly
accord to Leibniz's fictional infinitesimals the status of
mathematical entities that do not idealize anything in nature.
Bl{\aa}sj\"o's approach in \cite{Bl17} is largely compatible with such
a reading.

By contrast, the proponents of Ishiguroan interpretations deny the
Leibnizian infinitesimals the status of mathematical entities
violating Euclid V, Definition 4 (closely related to the so-called
Archimedean axiom).  Focusing on three examples, taken from the work
of Bassler (Section~\ref{s34}), Knobloch (Section~\ref{s51}) and
Arthur (Section~\ref{s5}), we signal some difficulties such approaches
encounter in interpreting Leibniz.

In Section~\ref{s34} we present the kind of problems which result from
conflating concepts related to the mathematical and the physical
realm.  In Section \ref{s51} we outline the difficulty which arises
from conflating, in mathematics, the Leibnizian notions of bounded and
unbounded infinite.  Finally, a situation where both kinds of
conflation are present is analyzed in Section~\ref{s5}.

\section{Creatures and infinitesimals}
\label{s34}

The locution \emph{syncategorematic infinitesimal} occurs nowhere in
Leibniz's known work.  It is coined on the model of the expression
\emph{syncategorematic infinite}, a crucial concept in Leibniz's
thought.

\subsection{Syncategorematic infinitesimals?}

Applying to the infinite the Mediaeval distinction between
categorematic and syncategorematic use of a term (see
Section~\ref{s6}), one can say that Leibniz's indefinite i.e.,
unbounded infinite (see \cite[Section\;4.3]{22a}) is syncategorematic,
because it does not refer to any properly infinite entity.  Rather, it
means that given a finite entity, such as a number or a segment, one
can always take a bigger finite one, and so on, indefinitely.  In
other words, the quantitative infinite must be intended always in a
\emph{distributive} way, so that for Leibniz there is no
\emph{collective} quantitative infinite.  On this definition see for
instance Bassler's \cite[p.\;855, note\;15]{Ba98}, singled out by
Richard Arthur for special praise:
\begin{enumerate}\item[]
``Bassler's erudite footnote on the syncategorematic and categorematic
  in his 1998, 855, n.\,15.'' (Arthur \cite{Ar15}, 2015, p.\,145,
  note\;15)
\end{enumerate}

On the relation between Leibniz's syncategorematic infinite and
Aristotle's potential infinite see \cite{Ug22}.

A similar notion of syncategorematicity is applied by Arthur, Bassler,
and others to infinitesimals.  Just as a syncategorematic infinite
refers to the procedure of indefinitely increasing a finite quantity,
which remains at all times finite, so also syncategorematic
infinitesimals refer to the procedure of indefinitely decreasing a
finite (more precisely, \emph{assignable}) quantity, which remains at
all times finite.  On this reading, infinitesimals are a shorthand for
Archimedean exhaustion-style procedures.

The thesis attributing such syncategorematic infinitesimals to Leibniz
depends on two main assumptions: Leibniz' alleged claims of the
non-existence of infinitesimals, and the identification of unbounded
infinities (\emph{infinita interminata}) as reciprocal of
infinitesimals.  In this section we discuss Bassler's version of the
first assumption, starting with a key passage from a 20\;june 1702
letter to Varignon, where Leibniz made the following intriguing
comment:

\begin{enumerate}\item[]
Je croy qu'il n'y a point de creature au dessous de la quelle il n'y
ait une infinit\'e de creatures, cependant je ne crois point qu'il y
en ait, ny m\^eme qu'il y en puisse avoir d'infiniment petites et
c'est ce que je crois pouvoir demonstrer.  \cite[p.\,110]{Le02c}
\end{enumerate}
The French noun \emph{cr\'eature} (spelled \emph{creature} by Leibniz)
is feminine.  The qualifier \emph{infiniment petites} accords with
\emph{creatures} and is therefore also feminine.  Thus, Leibniz is
asserting that there could not be infinitesimal (material) creatures.
A fragment of Leibniz's sentence is quoted by Bassler in translation.
We note the following four points.

\subsection{Infinitely small things}
\label{s21}

Bassler suppresses the beginning of the sentence and reproduces only
its conclusion as follows:
\begin{enumerate}\item[]
[Leibniz] believes he can prove `that there are not, nor could there
be, any infinitely small things'.  \cite[p.\;873]{Ba98}
\end{enumerate}
The sentence fragment quoted (in translation) by Bassler, taken out of
context, may suggest that Leibniz asserts the impossibility of
infinitesimals.  Read in context, it is clear that the sentence merely
asserts the impossibility of \emph{material} infinitesimals.%
%
%

Bassler is not alone in seeking to exploit this letter to Varignon in
support of a thesis denying the status of mathematical entities to
infinitely small quantities.  A similar truncation maneuver occurs in
Schubring, who purports to quote Leibniz as follows:
\begin{enumerate}\item[]
I do not believe that there are or even that there could be
\emph{infinitely small quantities}, and that is what I believe to be
able to prove. It is that simple substances (that is, those that are
not aggregated objects) are truly indivisible, but they are not
material and are but principles of action.  (Leibniz as rendered by
Schubring in \cite{Sc05}, 2005, p.\;173, lines\;-12
through\;\mbox{-9}; emphasis added)
\end{enumerate}
Notice that Schubring has Leibniz speak of ``infinitely small
quantities'' in this passage, whereas Leibniz himself only mentioned
``creatures infiniment petites.''  The modification is convenient for
scoring an anti-infinitesimal point, but distorts Leibniz's intention.
Additional instances of truncation of Leibnizian passages resulting in
distorting their meaning are presented in \cite{16a}.

\subsection{Bassler \emph{vs} Jesseph}

Bassler criticizes Jesseph (\cite{Je89}, 1989) in the following terms:
\begin{enumerate}\item[]
[Jesseph's] suggestion that we can use Leibniz's attitude concerning
the segregation of metaphysical from mathematical problems as a
resolution of the `ambiguity' in Leibniz's use of infinitesimals
{\ldots}\ seems to me hand-waving at best and circular at worst.
\cite[p.\;868, continuation of note\;57]{Ba98}
\end{enumerate}
According to Bassler, ``The problem, of course, is that such a
`radical thesis' itself requires a \emph{metaphysical} justification''
(ibid.; emphasis in the original).  Bassler provides no further
details concerning `the problem', or why a `segregation of
metaphysical from mathematical' would be so problematic in his
opinion.  

In mature Leibniz, one finds three distinct realms: (1) the real
metaphysical realm of substances or monads, (2) the phenomenal
physical realm; (3) the ideal mathematical realm.  The desegregation,
or more precisely conflation, of these realms is what led Bassler to
misrepresent the Leibnizian position, as noted in Section~\ref{s21}.

\subsection{Where do infinitesimals fit in?}
\label{s53b}

Bassler lodges the following claim against infinitesimals:
\begin{enumerate}\item[]
[E]xisting infinitesimals, whether ideal or real, fit nowhere.  This
is the ultimate consequence of Leibniz's understanding that in the
realm of quantity there are no infinite or infinitesimal
\emph{wholes}: the quantitative infinite is identified with the
indefinite or unlimited.  (\cite{Ba98}, p.\;872; emphasis on `wholes'
added)
\end{enumerate}
The claim suffers from a conflation of quantity and multitude.%
\footnote{A similar conflation in Rabouin and Arthur \cite{Ra20} is
  analyzed in \cite[Section\;4.2]{22a}.}
Leibniz held that infinite \emph{multitude} taken as a \emph{whole}
(i.e., \emph{infinita interminata}) contradicts the part-whole
principle (see~\cite[Section 4.1]{22a}).
However, Bassler provides no evidence that Leibniz viewed
infinitesimals and infinite \emph{quantities} (\emph{infinita
  terminata}) as contradictory.

\subsection{Relatively more freestanding symbols}
\label{s44}

Bassler proposes a syncategorematic reading of infinitesimals similar
to Ishiguro's:
\begin{enumerate}\item[]
Leibniz's recourse to these incomparable quantities is exclusively
symbolic: they are mathematical symbols whose true interpretation can
only be given in terms of their replacement, in any given context, by
other, relatively more freestanding symbols.  {\ldots} Consequently
Leibniz will later declare that infinitesimals and infinitely large
magnitudes are ``convenient fictions.''$^{62}$ \cite[p.\;869]{Ba98}
\end{enumerate}
Here Bassler seeks to interpret fictionality in terms of
syncategorematicity.  Bassler's footnote 62 attached to this passage
is significant for our purposes.  The footnote contains the following
passage from Leibniz's letter to Masson:
\begin{enumerate}\item[]
The infinitesimal calculus is useful with respect to the application
of mathematics to physics; however, that is not how I claim to account
for the nature of things [\emph{la nature des choses}].  For I
consider infinitesimal quantities to be useful fictions.  (Leibniz
\cite{Le16} as translated by Ariew in \cite{Le89a}, p.\;230)
\end{enumerate}
However, Bassler's interpretation of the letter to Masson involves a
conflation of mathematics and nature.  A similar conflation occurs in
the work of Rabouin and Arthur \cite{Ra20}; for details see
\cite[Section~2.2]{22a} as well as Section \ref{s69} below.


\section
{Knobloch's evolution on infinitesimals, bounded infinity}
\label{s51}

The definition of syncategorematic infinitesimals is grounded in the
assumption that for Leibniz the counterparts of infinitesimals are
unbounded infinities.  However, much evidence to the contrary exists
in Leibniz's oeuvre.  He specifies repeatedly that the counterparts of
infinitesimals are bounded infinities (\emph{infinita terminata}),
rather than unbounded infinities (\emph{infinita interminata}).

In this section we trace Knobloch's analysis of the bounded infinite
in Leibniz, and his evolving views concerning the Leibnizian
infinitesimal geometry and calculus.  

\subsection{Early writings}
\label{s31}

In Knobloch's early writings on Leibniz, infinitesimals were less than
any assignable quantity, Breger was criticized for claiming that
bounded infinities were not important later on, and there was no talk
of variables (and certainly not of Weierstrass and the
\emph{Epsilontik}).

\medskip
\textbf{1989.}  Thus, in 1989, Knobloch stresses the importance of the
distinction between bounded and unbounded infinite in the following
terms:
\begin{enumerate}\item[]
Leibniz discute un probl\`eme math\'ematique qui est en m\^eme temps
un probl\`eme philosophique: l'infini dans les math\'ematiques, en
particulier la diff\'erence entre l'infini et l'intermin\'e {\ldots}
%
%
Cette diff\'erence est extr\^emement importante.  L'infini est une
quantit\'e qui est plus grande qu'une quantit\'e quelconque assignable
ou qu'une quantit\'e qu'on peut d\'esigner par des nombres.
\cite[p.\;166]{Kn89}
\end{enumerate}
Here \emph{l'infini} is the \emph{infinitum terminatum}, defined as
being greater than any assignable quantity.%
\footnote{\label{f60}Here Knobloch offers no reservations concerning
  such a Leibnizian definition of \emph{infinitum terminatum} in terms
  of assignable quantities; cf.\;note~\ref{f61}.}
Meanwhile, \emph{l'intermin\'e} is the unbounded infinity (see
\cite[Section~2.2]{21a}).
%
%
Significantly, there is no mention of variables in the 1989 text.


\medskip
\textbf{1990.}  Knobloch analyzes Leibniz's manuscript \emph{De
  Quadratura Arithmetica} (1676) in a 1990 text, and writes:
\begin{enumerate}
\item
L'action de l'esprit par laquelle nous mesurons les espaces infinis
s'appuie sur une fiction, sur une ligne, qui est \emph{termin\'ee}
selon la pr\'esupposition, mais infinie.  \cite[p.\;39]{Kn90}
(emphasis added)
\item
Si l'abscisse~$\mu(\mu)$ est infiniment petite, l'ordonn\'ee
$(\mu)\lambda$ est infiniment longue, c'est-\`a-dire plus grande
qu'une droite quelconque d\'esignable, le rectangle effectu\'e par la
droite infinie et la droite infiniment petite est \'egal \`a un
carr\'e fini constant selon la nature de l'hyperbole.
\cite[p.\;41]{Kn90}
\item
{}[Leibniz] dit dans une note marginale d'un \'ecrit de Schuller: J'ai
fait toujours une diff\'erence entre \emph{immense} et intermin\'e,
entre celui auquel on ne peut ajouter rien et celui qui est plus grand
qu'un nombre assignable.  (\cite[p.\;42]{Kn90}; emphasis added)
\item
L'infiniment petit est d\'efini par `plus petit qu'une grandeur
assignable' (\emph{minus quavis assignabili quantitate}).%
\footnote{\label{f61}This definition of infinitesimal as smaller than
  every assignable quantity is, appropriately, not accompanied by any
  vanishing claims, unlike its later occurrences in Knobloch's
  writing; see main text at note~\ref{f65}.}
\cite[p.\;49]{Kn90}
\end{enumerate}
Here the \emph{immense}, referring to \emph{infinita terminata}, is
again contrasted with unbounded infinity.  Knobloch's 1990 article
similarly contains no description of infinitesimals either as
variables or as logical fictions concealing alternating quantifiers.

\subsection{Knobloch \emph{vs} Breger}
\label{s72}

The issue of the bounded infinite was debated by Breger and Knobloch
in the 1990s.

\medskip
\textbf{1994.}  In 1994, Knobloch takes issue with Breger in defense
of the bounded infinite:
\begin{enumerate}\item[]
All the more surprising is the fact that a contrary remark is found in
the secondary literature which appeared very recently (Breger 1990,
65): `Leibniz seems to have called in question the utility of the
bounded infinite already during his stay in Paris.  Later on he seems
to have abandoned the bounded infinite in mathematics.'  The contrary
is correct.  \cite[p.\;268]{Kn94}
\end{enumerate}
The reference is to Breger's article (\cite{Br90}, 1990).  Knobloch
then goes on to provide evidence refuting Breger's claims.  Knobloch
writes further:
\begin{enumerate}\item[]
[Leibniz] accepts different orders of infinite lines and different
infinitely small lines.  If we choose an arbitrary finite line
$\ell_f$, we find an infinitely small line~$\ell_{is}$ and an infinite
line~$\ell_i$, so that we obtain the equation {\ldots}
$\ell_f^2=\ell_i \cdot \ell_{is}$   \cite[p.\;268]{Kn94}
\end{enumerate}
The passage appears to distinguish between three orders of magnitude
among Leibnizian \emph{lines} (i.e., segments) -- infinitesimal,
finite, and infinite -- typical of non-Archimedean systems.  The 1994
article similarly contains no description of infinitesimals as
variables.

\medskip
\textbf{1999.}  Knobloch has been writing about Leibnizian
infinitesimal geometry and analysis at least since (\cite{Kn78},
1978).  Over two decades later in 1999, one finds Knobloch for the
first time describing infinitesimals as variable quantities:
\begin{enumerate}\item[]
Leibniz defined `lines', that is linear indivisibles as infinitely
small rectangles, that is as variable quantities:%
\footnote{\label{f62}Similar claims concerning variable quantities
  appear in a 2018 article by Knobloch; see Section~\ref{s67}.}
His starting point was a quantification of the notion of indivisibles.
\cite[p.\;95]{Kn99}
\end{enumerate}
Yet in the same paper one still finds the reading of ``not appearing
in nature'' type:
\begin{enumerate}\item[]
The demonstrations of the forty-five theorems that follow [in DQA] are
based on this concept of infinitely small and infinite quantities,
that is, quantities which are, according to his definition, positive,
but smaller than any given quantity or larger than any given quantity.
The decisive aspect is that they are quantities, fictive ones
certainly, because they were introduced by a fiction, but quantities
nevertheless.  \emph{It does not matter whether they appear in nature
or not}, because they allow abbreviations for speaking, for thinking,
for discovering, and for proving.  (ibid.; emphasis added)
\end{enumerate}
Knobloch's 1999 interpretation of fictionality in terms of not
necessarily ``appearing in nature'' is consonant with the reading
developed in \cite{13f}, \cite{21a}, and \cite{22a}.

\medskip
\textbf{2002.} In 2002, Knobloch still quotes the Leibnizian
\begin{enumerate}\item[]
`law of continuity': The rules of the finite remain valid in the
domain of the infinite.  \cite[p.\;67]{Kn02}
\end{enumerate}
To escape vacuity, such a law would apparently need to refer to a
notion of an infinite beyond merely a variable finite---for why
wouldn't the rules of the finite naturally remain valid in the domain
of the variable finite, even without invoking any special laws or
``metaphysical principles of continuity'' (De Risi \cite{De19}, 2019,
p.\;124)?

\subsection{Quantifiers, Weierstrass, and the \emph{Epsilontik}}
\label{s67}

In Knobloch's later writings, one detects a change of attitude toward
Leibnizian infinitesimals.

\medskip
\textbf{2008.}  Thus, by 2008, alternating quantifiers make their
appearance, as Knobloch compares Leibniz and Euler:

\begin{enumerate}\item[]
We might sum up the crucial difference between Euler's and Leibniz's
definitions of infinitely small quantities in the following way: Let
$i$ be an infinitely small quantity,~$gq$ a given quantity,~$aq$ an
assignable quantity.
\medskip
\begin{enumerate}\item[]
Leibniz: For all~$gq > 0$ there is an~$i (gq) > 0$ so
that~\mbox{$i(gq)<gq$} \\~$\Rightarrow i (gq)$ is a variable quantity.

\medskip\noindent
Euler: For all~$i$ and for all~$aq > 0: \quad i < aq$ \\~$\Rightarrow
i = 0$.
\end{enumerate}
(Knobloch \cite{Kn08}, 2008, p.\;180)
\end{enumerate}
Knobloch's analysis of Leibnizian infinitesimals in terms of the
quantifier phrase ``For all~$gq > 0$ there is an~$i (gq) > 0$ so that
\mbox{$i (gq) < gq$}'' echoes Ishiguro's syncategorematic alternating
quantifier reading; see Section \ref{s6}.  Knobloch's claim of such a
dramatic discontinuity in attitude toward infinitesimals between
Leibniz and Euler is dubious.%
\footnote{Cauchy has been similarly the subject of dubious
  interpretations aimed at denying that he worked with \emph{bona
    fide} infinitesimals; for details see \cite{19a}, \cite{20a},
  \cite{20b}.}
A case in favor of Leibniz--Euler continuity is argued in Bair et
al.\;(\cite{17b}, 2017).

\medskip
\textbf{2012.}  In later Knobloch, both Weierstrass and the
\emph{Epsilontik} make their appearance, accompanied by (1) a denial
of inassignables, and (2)~insertion of alternating quantifiers.  Thus,
in 2012 we learn that
\begin{enumerate}\item[]
From 1673 he tried different possibilities like smallest,
unassign\-able magnitude, smaller than any assignable quantity.  He
rightly rejected all of them because there are no smallest quantities
and because a quantity that is smaller than any assignable quantity is
equal to zero or nothing.%
\footnote{\label{f63}A quantity~$Q$ smaller than any assignable
  quantity is of course necessarily zero provided~$Q$ itself is
  assignable.  Such an observation is the last step of a typical
  argument by exhaustion.  However, Knobloch's conclusion does not
  follow if~$Q$ is not assumed to be assignable.  This enables
  inassignable infinitesimals to be both smaller than every assignable
  quantity \emph{and} nonzero.  See \cite[Section~4.5]{22a} for a
  geometric illustration in terms of hornangles.  For an analysis of a
  related misconception in Rabouin \cite{Ra15} see
  \cite[Section~3.2]{22a}.}
{\ldots} Now [Leibniz] had to answer the question: What does it mean
to be infinitely small? Still in 1673 he gave an excellent answer:
infinitely small means smaller than any given quantity.  (Knobloch
\cite{Kn12}, 2012, pp.\;20--21)
\end{enumerate}
Knobloch elaborates:
\begin{enumerate}\item[]
He again used the mode of possibility and introduced a consistent
notion.  Its \emph{if-then structure} -- if somebody proposes a
quantity, then there will be a smaller quantity -- rightly reminds the
modern reader of Weierstra\ss's~$\epsilon$-$\delta$-language.
Leibniz's language can be translated into Weierstra\ss's language.
(Knobloch \cite{Kn12}, 2012, p.\;21; emphasis added)
\end{enumerate}
The postulation of such an ``if-then structure'' is akin to Ishiguro's
interpretation.  Knobloch's attribution of such a reading to Leibniz
as early as 1673 clashes with Arthur's claim that the switch to
syncategorematics occurred between february and april of 1676.%
\footnote{\label{f64}Thus, Arthur writes: ``For some time, Leibniz
  appears to have hesitated over this interpretation, and as late as
  February 1676 he was still deliberating about whether the success of
  the hypothesis of infinities and the infinitely small in geometry
  spoke to their existence in physical reality too. But by April
  [1676], the syncategorematic interpretation is firmly in place''
  \cite[p.\;559]{Ar13}.}
For additional details see \cite[Section~3.1]{22a}.

\medskip
\textbf{2018.}  In 2018, Knobloch writes: ``Such a difference or such
a quantity [i.e., an infinitesimal] necessarily is a variable
quantity'' \cite[p.\;12]{Kn18}.  In the same text Knobloch notes that
\begin{enumerate}
\item[]
``the error is smaller than any assignable error and therefore zero
  {\ldots} Such an error necessarily is equal to zero as Leibniz
  rightly states. For if we assume that such an error is unequal to
  zero it would have a certain value.  But this implies a
  contradiction against the postulate that the error has to be smaller
  than any assignable quantity, that is, also smaller than this
  certain value.  Yet, Leibniz explicitly calls such errors infinitely
  small: \emph{We should not try to make things seem better}'' (ibid.;
  emphasis added).
\end{enumerate}
The italicized comment appears to imply that something could be better
about the Leibnizian remarks on infinitesimals.  Actually things seem
rather good to the present authors.  What Knobloch is bothered by is
an apparent tension between the fact that the error is supposed to
turn out to be zero rather than merely infinitely small.  But the
desired conclusion will be correct if the quantity being estimated is
assumed to be assignable.%
\footnote{This point is explained in note~\ref{f63}.}

Also in 2018, we learn that
\begin{enumerate}\item[]
[Leibniz] himself noticed that the definition ``smaller than any
assignable'' is useless because such a quantity must be necessarily
equal to zero.%
\footnote{\label{f65}The misconception concerning infinitesimals being
  necessarily zero has already been analyzed in notes~\ref{f62} and
  \ref{f63}.  Note that in 1989 and 1990 Knobloch presented such
  definitions unaccompanied by any vanishing claims; see
  notes~\ref{f60} and~\ref{f61}.}
\cite[p.\;54]{Kn18b}
\end{enumerate}
But did Leibniz ever notice such a thing?  The impossibility of a
nonzero infinitesimal would apparently imply also the impossibility of
its reciprocal, the \emph{infinita terminata}.  Yet \emph{infinita
  terminata} was staunchly defended by early Knobloch in 1994 against
Breger (see Section~\ref{s72}).

\medskip
\textbf{2022.}  Most recently (2022?), we learn that
\begin{enumerate}\item[]
An infinitely small quantity is a variable quantity and can be
described in terms of the Weierstrassian epsilon-delta language.
\cite{Kn22}
\end{enumerate}
It emerges that according to later Knobloch, the Leibnizian infinite
is not so \emph{immense} after all.  We see that the task of declaring
Leibniz a Weierstrassian and practitioner of the \emph{Epsilontik}
took Knobloch several decades of work.  However, his early writings
sampled in Sections~\ref{s31} and~\ref{s72} arguably reflect a more
accurate picture of the Leibnizian infinitesimal geometry and
calculus.

\section
{An Arthurian Cantor--Leibniz s\'eance}
\label{s5}

In a 2019 medium, Richard Arthur conjures up a 33-page Cantor--Leibniz
dialogue.  Here Arthur endeavors
\begin{enumerate}\item[]
[to] use this dialogue to argue that the Leibnizian actual infinite
constitutes a perfectly clear and consistent third alternative in the
foundations of mathematics to the usual dichotomy between the
potential infinite (Aristotelianism, intuitionism) and the transfinite
(Cantor, set theory), and one that avoids the paradoxes of the
infinite.  \cite[p.\;72]{Ar19b}
\end{enumerate}
Concerning the actual infinite, Arthur claims in 2019 that
\begin{enumerate}\item[]
[A]ctual infinite syncategorematically understood {\ldots} as a
distributive whole \ldots\;constitutes a perfectly clear and
consistent third alternative in the foundations of mathematics.
\cite[p.\;106]{Ar19b}
\end{enumerate}
Arthur's proposal of such alternative foundations is based on a
passage from a supplement to Leibniz's letter to des Bosses.  Look and
Rutherford translate this letter to des Bosses (number 10 in their
collection~\cite{Le07}) as well as the ``supplementary study'' which
contains the passage about ``actual infinite in the sense of a
distributive whole.''  In footnote\;8 on page 409, Look and Rutherford
note that this passage was ``crossed through by Leibniz.''  This
crossed-out comment is what Arthur is basing his ``perfectly clear
alternative foundations for mathematics'' on.  The reviewer for
Zentralblatt observes:
\begin{enumerate}\item[]
[T]wo concerns of this paper seems to me not correct: 1) the idea that
through the concept of the actually infinite division of matter and
through the highlighted meaning of syncategorematic infinite a
consistent conception of the actual infinite may be obtained; 2)
Cantor's conception and theory of transfinite is inconsistent.
(Bussotti \cite{Bu19})
\end{enumerate}

\subsection{Evolution of Arthur's views}

Arthur's views concerning mathematical actual infinity in Leibniz
appear to have undergone a significant evolution.  In 2009, he claimed
the following:
\begin{enumerate}\item[]
In many quarters, to boot, [Leibniz's] mature view is seen as somewhat
unfortunate, especially since the work of Abraham Robinson and others
in recent years, which has succeeded in rehabilitating infinitesimals
as \emph{actual}, non-Archimedean entities.  A dissenting view has
been given by Ishiguro, who argues (in my opinion, persuasively) that
Leibniz’s syncategorematic interpretation of infinitesimals as
fictions is a conceptually rich, consistent, finitist theory, well
motivated within his philosophy, and no mere last-ditch attempt to
safeguard his theory from foundational criticism.  \cite[p.\;11]{Ar09}
(emphasis on `actual' added)
\end{enumerate}
Accordingly, it is only in Robinson-influenced quarters that there is
any talk of infinite actuality; in Ishiguro's ``persuasive'' view,
infinitesimals are part of a ``finitist theory.''  Similarly in 2013,
Arthur writes:
\begin{enumerate}\item[]
I take the position here (following Ishiguro 1990) that the idea that
Leibniz was committed to infinitesimals as \emph{actually infinitely
  small entities is a misreading}: his mature interpretation of the
calculus was fully in accord with the Archimedean Axiom.  {\ldots}
Infinitesimals {\ldots} really stand for variable finite quantities
that can be taken as small as desired.  \cite[p.\;554]{Ar13} (emphasis
added)
\end{enumerate}
Here again, actuality is reported to be a ``misreading.''

Starting no later than 2014, we begin to find references to the
syncategorematic actual infinite in mathematics:
\begin{enumerate}\item[]
There is an \emph{actual} infinite, but it must be understood
\emph{syncategorematically}: to say that a series has infinitely many
terms means not that it has~$N$ terms where~$N$ is greater than any
finite number, but that for any finite~$N$ it has more terms than
that, ‘more terms than one can specify’.  \cite[p.\;88]{Ar14}
(emphasis on `syncategorematically' in the original; emphasis on
`actual' added)
\end{enumerate}
Arthur adequately explains the point about `more terms than~$N$' but
provides no evidence for his \emph{actual} claim.  \emph{Pace} Arthur,
Robinson specifically emphasized that Leibniz accepted only potential
syncategorematic infinity of mathematical series:
\begin{enumerate}\item[]
[A]lthough Leibniz wishes to make it clear that he can dispense with
infinitely small and infinitely large numbers and at any rate does not
ascribe to them any absolute reality, he accepts the idea of potential
(‘syncategorematic’) infinity and finds it in the totality of terms of
a geometrical progression.  (Robinson \cite{Ro66}, 1966, p.\;263)
\end{enumerate}
A similar view is expressed in (Robinson \cite{Ro67}, 1967, p.\;35).%
\footnote{For a formalisation of the Leibnizian calculus in Robinson's
  framework, see \cite{17f}.}
Robinson's position is compatible with those of other Leibniz scholars
including Antognazza \cite{An15}; see \cite[Section~2.3]{22a}.

Arthur's claim concerning actual infinity in Leibniz involves a
conflation of the two realms of (1) mathematics and (2) nature, since
Leibnizian actual infinity (understood distributively) pertains only
to the latter.

\subsection{`Testament to the necessity of infinite number'}
\label{s69}

In 39 pages, Arthur's 2019 text firmly establishes that Leibniz would
have rejected Cantorian infinite cardinalities because they contradict
the part-whole principle.%
\footnote{A point already made by Ishiguro in 1990 (``Leibniz did not
  think that there should be what we call the cardinality of the set
  of all things'' \cite[p.\;80]{Is90}) as well as earlier authors; see
  e.g., Russell (\cite{Ru03}, 1903, p.\;210, end of item\;140).
  Russell claims further that ``if there were no infinite wholes, the
  word \emph{Universe} would be wholly destitute of meaning'' (ibid.).
  Russell accepted, as a matter of ontological certainty, Cantor's
  position both concerning the existence of infinite wholes and the
  non-existence of infinitesimals.  For an analysis of a dissenting
  opinion by contemporary neo-Kantians see \cite{13h}.}
This ``transcription'' \cite[p.\;72]{Ar19b} of an afterlife
Cantor--Leibniz dialogue seeks to graft a passage from the letter to
Masson upon the infinitesimal calculus.  A fragment of the dialogue
runs as follows:

\begin{enumerate}\item[]
\textbf{Cantor}: {\ldots} And what about your Infinitesimal Calculus?
You must agree that this is a testament to the necessity of infinite
number for mathematics!

\vspace{3pt}
\noindent 
\textbf{Leibniz}: On the contrary, in spite of my Infinitesimal
Calculus, I admit no genuine (v\'eritable) infinite number, even
though I confess that the multiplicity (multitude) of things surpasses
every finite number {\ldots} 

\vspace{3pt}
\noindent
(Cantor and Leibniz as transcribed by Arthur in \cite{Ar19b}, 2019,
p.\;74)
\end{enumerate}
In his note 11 on page 74, Arthur sources the Leibnizian passage above
in the 1716 letter to Masson.  However, the Leibnizian passage is
taken out of context.  Arthur exploits the passage to suggest that
Leibniz abstained from using infinite quantities in his calculus.  The
passage, taken out of context, can be read that way.  In fact Leibniz
does not mean to suggest such abstinence.  While we don't disagree
with Arthur that Leibniz may have rejected Cantorian infinite
\emph{cardinalities}, he did use infinite \emph{quantities} that are
\emph{mentis fictiones} (that may be impossible in nature), as argued
in \cite{22a}.  Thus Arthur's ``pastiche'' \cite[Note\;1]{Ar19b} tends
to misrepresent Leibniz' philosophical stance with regard to
infinities.

Leaving aside the question whether Georg Cantor -- who referred to
non-Archi\-me\-dean number systems as the ``infinite cholera
bacillus''%
\footnote{In the original German: \emph{diesen infinit\"aren
Cholera-Bacillus} \cite[p.\;505]{Me}.}
of mathematics%
\footnote{\label{f68}Cantor names Thomae specifically.  See Ehrlich
  (\cite{Eh06}, 2006, p.\;55) for an analysis of Thomae's
  non-Archimedean mathematics; see also B{\l}aszczyk and Fila
  \cite{Bl20}.  In the same letter, Cantor writes:
\begin{enumerate}\item[]
Professor Veronese (who, if I am not mistaken, emerged from the famous
school of Professor Felix Klein, whom he also revered as a mighty
mathematical hero), on the shoulders of those German predecessors in
connection with the superficial and unclear ideas of his just
mentioned teacher, finally built up his new theory of the actual
infinitely small segments and claimed to have created by these a new
foundation for geometry. (\cite[p.\;505]{Me}; translation ours)
\end{enumerate}
Cantor's disdainful comments on Klein are all the more surprising
since Klein was the one who came to Cantor's rescue when the latter
was unable to publish his set-theoretic research in Berlin, and
rapidly published a series of papers by Cantor in \emph{Mathematische
  Annalen}.  For further details on Klein's modernity see \cite{18b}.}
-- would be likely to press their case (i.e., the case of the infinite
numbers in the calculus) in a s\'eance after his demise,%
\footnote{Cantor was, of course, interested in infinite cardinalities
  and sometimes referred to them as numbers, but these are clearly not
  the infinite numbers of Infinitesimal Calculus that he is presented
  as a champion of in this fragment of Arthur's s\'eance.  Leibniz
  already clearly understood the difference between multitude and
  magnitude.}
one notes a conflation of the issue of Leibnizian opposition to
infinite wholes (entailing, as argued by Arthur, a rejection of
Cantorian infinite cardinalities); and the well-founded fictionality
of his infinitesimal and infinite quantities.

The unique form of Arthur's medium, and the artistic licence it grants
itself, do not yet justify playing loose with the philosophical
positions of both Cantor and Leibniz.

\section{Conclusion}

We presented three case studies illustrating the fact that current
received Leibniz scholarship routinely denies the Leibnizian
infinitesimals the status of mathematical entities violating Euclid V
Definition~4.  We highlighted some difficulties that such
interpretations encounter relative to the Leibnizian texts.  Such
difficulties do not arise for interpretations of Leibnizian
infinitesimals that accord them the status of mathematical entities,
that have recently appeared in the scholarly literature; see e.g.,
\cite{21a}, \cite[pp.\;620, 641]{Es21}, \cite{22a}, \cite{22b}.

\section*{Acknowledgment}

The influence of Hilton Kramer (1928--2012) is obvious.

\end{document}